\documentclass[a4paper,11pt]{amsart}
\usepackage{amsmath,amsthm,amssymb}
\usepackage[arrow,matrix]{xy}
\usepackage{pstricks,pst-grad, pst-node}

\theoremstyle{plain}
\newtheorem{theorem}{Theorem}[section]
\newtheorem{proposition}[theorem]{Proposition}

\newtheorem{problem}[theorem]{Problem}

\theoremstyle{definition}

\newtheorem{definition}[theorem]{Definition}

\newtheorem{example}[theorem]{Example}

\newcommand{\Q}{\mathbb{Q}}

\newcommand{\CP}{\mathbb{C}P}

\newcommand{\cZ}{\mathcal{Z}}
\newcommand{\cP}{\mathcal{P}}
\newcommand{\cQ}{\mathcal{Q}}

\DeclareMathOperator{\Tor}{Tor}

\DeclareMathOperator{\bideg}{bideg}
\def\k{\mathbf{k}}

\begin{document}%
\title[Moment-angle manifolds with the same bigraded Betti numbers]{Different moment-angle manifolds arising from two polytopes having the same bigraded Betti numbers}
\author{Suyoung Choi}
\address{Department of Mathematics, Ajou University,
San 5, Woncheondong, Yeongtonggu, Suwon 443-749, Korea}
\email{schoi@ajou.ac.kr}
\thanks{The author was partially supported by the Basic Science Research Program through the National Research Foundation of Korea (NRF) funded by the Ministry of Education, Science and Technology (2011-0024975).}

\date{\today}
\begin{abstract}
    Two simple polytopes of dimension $3$ having the identical bigraded Betti numbers but non-isomorphic Tor-algebras are presented. These polytopes provide two homotopically different moment-angle manifolds having the same bigraded Betti numbers. These two simple polytopes are the first examples of polytopes that are (toric) cohomologically rigid but not combinatorially rigid.
\end{abstract}

\maketitle

\tableofcontents

\section{Introduction}%

A convex polytope of dimension $n$ is called \emph{simple} if there are exactly $n$ facets (codimension-one face) meeting at each vertex.
Let $P$ be an $n$-dimensional simple convex polytope with $m$ facets $F_1, \ldots, F_m$. Consider an $m$-dimensional real compact torus $T^m$, and denote the $i$-th coordinate subgroup of $T^m$ by $T_i$.

\begin{definition}
  Consider the following equivalence relation on $T^m \times P$:
  $$
    (t, p)\sim (t', p') \Longleftrightarrow p=p', t't^{-1} \in \bigoplus_{F_i \ni p} T_i.
  $$
  Then, the quotient space
  $$
    \cZ_P = (T^m \times P) / \sim
  $$ is called the \emph{moment-angle manifold} of $P$ and is denoted by $\cZ_P$.
\end{definition}

It is noted that $\cZ_P$ is indeed a manifold of dimension $m+n$ (see \cite[Lemma 6.2]{BP}), and the formula $s \cdot [t,p] = [st, p]$ defines a natural $T^m$-action on $\cZ_P$ with orbit space $P$. The moment-angle manifold was introduced in \cite{DJ} as a space that has the following universal property: for every quasitoric manifold (the definition will be given below) $\pi : M \to P$, there is a principal $T^{m-n}$-bundle $\cZ_P \to M$ whose composite map with $\pi$ is the orbit map $\cZ_P \to P$. Hence, it is one of the key concepts in toric topology, and it is very important to study the topology of $\cZ_P$.

A formula for the cohomology of $\cZ_P$ has already been established. Let $\k$ be a field. The Tor-algebra of $P$, denoted by $\Tor^{\ast,\ast}_A(\k(P),\k)$, is a finite-dimensional \emph{bigraded} $\k$-algebra. (The explicit definition will be given in Section~\ref{sec:tor-algebra}). Note that the cohomology algebra $H^\ast(\cZ_P,\k)$ of $\cZ_P$ inherits a canonical bigrading from the Eilenberg-Moore spectral sequence for the fibration
$$
\xymatrix{
   \cZ_K  \ar[r] \ar[d]^{=} &    \cZ_K  \ar[r] \ar[d] & \ast \ar[d]  \\
   \cZ_K \ar[r] & ET^m \times_{T^m}    \cZ_K \ar[r] & BT^m,}
$$ where $ET^m$ is a contractible space on which $T^m$ acts freely, and $BT^m = ET^m/T^m$.
Buchstaber and Panov (Theorems 7.6 and 7.7 in \cite{BP}) showed that $H^{\ast,\ast}(\cZ_P,\k)$ and $\Tor^{\ast,\ast}_A(\k(P),\k)$ are isomorphic as bigraded $\k$-algebra.

In this study, it is assumed that $\k$ is the field of rational numbers $\Q$. Let $\beta^{-i,2j}(P)$ denote the bigraded Betti numbers of the Tor-algebra of $P$ (simply, the bigraded Betti numbers of $P$), that is,
$$\beta^{-i,2j}(P) = \dim_\Q \Tor^{-i,2j}_A(\Q(P),\Q).$$
It is to be noted that the bigrading structure should have more information than the usual (mono)grading structure. Hence, it is natural to ask how much information on the topology of $\cZ_P$ the bigraded Betti numbers have. Actually, in all known examples (before this paper) of combinatorially different polytopes with the same bigraded Betti numbers (such as vertex truncations of simplices), the moment-angle manifolds are also diffeomorphic.
It should also be noted that the bigraded Betti numbers of $H^{\ast,\ast}(\cZ_P,\Q)$ are not necessary for them to be  topological invariants, although the usual Betti numbers $\beta^{p}=\sum_{p=-i+2j} \beta^{-i,2j}(P)$ are topological invariants.

From this viewpoint, Panov presented the following problem at the conference on toric topology held in Osaka in November 2011.

\begin{problem}
Let $P$ and $Q$ be two simple polytopes. Is it true that
$$\cZ_P \cong \cZ_Q \Longleftrightarrow \beta^{-i,2j}(P) = \beta^{-i,2j}(Q) \quad \text{ for all $i,j$?}$$
Here, $\cong$ may mean ``homotopy equivalent,'' ``homeomorphic,'' or ``diffeomorphic''.
\end{problem}

In this paper, we answer the ``if'' part of the problem negatively for all categories, namely, there exist two simple polytopes (say $\cP$ and $\cQ$) with the same bigraded Betti numbers, satisfying $H^\ast(\cZ_\cP) \not\cong H^\ast(\cZ_\cQ)$ as rings. Such polytopes are shown in Figure~\ref{fig:P_and_Q}. Note that both polytopes are $3$-dimensional simple polytopes having $11$ facets.
\begin{figure}\label{fig:P_and_Q}
\begin{pspicture}(0,0)(5,5)
\pspolygon(2.5,0.2)(0.5,1.1)(0.2,1.9)(0.5,3.2)(1.4,4.4)(2.5,4.8)(3.6,4.4)(4.5,3.2)(4.8,1.9)(4.5,1.1)
\psline(0.5,3.2)(2.5,2.5)\psline(4.5,3.2)(2.5,2.5)
\psline(0.2,1.9)(1.9,1.4)\psline(4.8,1.9)(3.1,1.4)
\psline(2.5,2.5)(2.5,2)\psline(2.5,0.2)(2.5,1)
\pspolygon(2.5,1)(1.9,1.4)(2.5,2)(3.1,1.4)
\psline[linestyle=dotted] (1.4,4.4)(1.4,3.1)
\psline[linestyle=dotted] (2.5,4.8)(2.5,3.5)
\psline[linestyle=dotted] (3.6,4.4)(3.6,3.1)
\psline[linestyle=dotted] (0.5,1.1)(1.4,3.1)(2.5,3.5)(3.6,3.1)(4.5,1.1)
\end{pspicture}
\begin{pspicture}(0,0)(5,5)
\pspolygon(3.6,0.4)(4.7,1.2)(4.9,2.5)(4.4,3.9)(3.5,4.7)(1.5,4.7)(0.6,3.9)(0.1,2.5)(0.3,1.2)(1.4,0.4)
\psline(4.4,3.9)(3.4,3) \psline(0.6,3.9)(1.6,3)
\psline(4.9,2.5)(3.6,1.6) \psline(0.1,2.5)(1.4,1.6)
\psline(3.6,0.4)(3.6,1.6)(3.1,2.2) \psline(1.4,0.4)(1.4,1.6)(1.9,2.2)
\pspolygon(3.1,2.2)(3.4,3)(1.6,3)(1.9,2.2)
\psline[linestyle=dotted] (4.7,1.2)(3.5,2.7)
\psline[linestyle=dotted] (0.3,1.2)(1.5,2.7)
\psline[linestyle=dotted] (3.5,4.7)(3.5,2.7)(1.5,2.7)(1.5,4.7)
\end{pspicture}
 \caption{$\cP$ and $\cQ$}
\end{figure}

Using the algebra program \emph{Macaulay2}\footnote{\emph{Macaulay2} can be downloaded from http://www.math.uiuc.edu/Macaulay2/}, one can see that they have the same bigraded Betti numbers. Table~\ref{Table:Betti_numbers_of_P_and_Q} presents the complete list of the bigraded Betti numbers of $\cP$ and $\cQ$.

\begin{table}
  \centering
\begin{tabular}{|c|ccccccccc|}
  \hline
  $j \setminus i$ & 0 & 1 & 2 & 3 & 4 & 5 & 6 & 7 & 8\\ \hline
  0 & 1&  &  &  &  &  &  & &  \\
  1 &  &  &  &  &  &  &  & & \\
  2 &  &28  &  &  &  &  &  & &\\
  3 &  &0  &105  &  &  &  &  & & \\
  4 &  &  &4  &166  &  &  &  & & \\
  5 &  &  &  &39  &123  &  &  & & \\
  6 &  &  &  &  &123  &39  &  & & \\
  7 &  &  &  &  &  &166  &4  & & \\
  8 &  &  &  &  &  &  &105  &0 & \\
  9 &  &  &  &  &  &  &  &28 &\\
  10&  &  &  &  &  &  &  & &\\
  11&  &  &  &  &  &  &  && 1\\
  \hline
\end{tabular}
  \caption{bigraded Betti numbers of $\cP$ and $\cQ$}
  \label{Table:Betti_numbers_of_P_and_Q}
\end{table}

\begin{theorem}\label{Theorem:Tor_is_not_decided_by_Betti_numbers}
The Tor-algebras of $\cP$ and $\cQ$ are not isomorphic as algebras. Furthermore, $H^\ast(\cZ_P)$ and $H^\ast(\cZ_Q)$ are not isomorphic as rings.
\end{theorem}

As an immediate corollary, it follows that the bigraded Betti numbers of the simple polytopes do not decide the homotopy type of the corresponding moment-angle manifold.

As by-products, the polytopes $\cP$ and $\cQ$ are one important examples in the toric rigidity problem for  simple polytopes as follows. A \emph{quasitoric manifold} is a closed smooth manifold of dimension $2n$ that admits a locally standard half-dimensional torus action $T^n$ whose orbit space is a simple polytope (see \cite{BP} and \cite{DJ}). A typical example of a quasitoric manifold is a complex projective space $\CP^n$ of complex dimension $n$ with the standard $T^n$-action whose orbit space is the $n$-simplex $\Delta^n$. Although the topology of a quasitoric manifold does not generally determine the combinatorial type of its orbit space, it sometimes does; for instance, only the $n$-simplex can be the orbit space of a locally standard $T^n$-action defined on $\CP^n$. Furthermore, since the cohomology ring $H^\ast(M)$ of a quasitoric manifold $M$ can be obtained from the face ring of its orbit polytope $P$, the relationship between the combinatorial type of $P$ and $H^\ast(M)$ is well established (see \cite{ch-pa-su10} and \cite{DJ}).

A simple polytope is said to be \emph{cohomologically rigid} if its combinatorial structure is decided by the cohomology ring of a supporting quasitoric manifold. A simple polytope is said to be \emph{combinatorially rigid} if its combinatorial structure is decided by the bigraded Betti numbers. By \cite[Proposition 3.8]{ch-pa-su10}, the bigraded Betti numbers of a simple polytope are determined by the cohomology ring of a supporting quasitoric manifold. Hence, any combinatorially rigid polytope (that supports a quasitoric manifold) is cohomologically rigid. However, the question of whether the converse holds has been open (see \cite[Section 6]{ch-ma-su11} for details).

\begin{problem}[Problem 6.6 in \cite{ch-ma-su11}]
  Find a polytope which is rigid cohomologically but not combinatorially in the set of simple polytopes.
\end{problem}

Here, we provide an answer to this problem.

\begin{theorem}\label{Theorem:cohom_does_not_imply_combi}
The polytopes $\cP$ and $\cQ$ are cohomologically rigid, but not combinatorially rigid.
\end{theorem}

\section{Tor-algebra of a simple polytope} \label{sec:tor-algebra}
we briefly review the definitions here, following \cite{BP}, where the reader may find additional details of the Tor-algebra of a simple polytope, and we present the properties of the multiplicative structure of the Tor-algebra, which are relevant to Section~\ref{sec:proof_of_main_theorem}.

Let $\k$ be a field, and let $A = \k[v_1, \ldots, v_m]$ be a finitely generated commutative graded algebra over $\k$. Then, $\k$ itself is an $A$-module via the map $A \to \k$ that sends each $v_i$ to $0$. Let $\Lambda[u_1, \ldots, u_m]$ denote an exterior algebra on $m$ generators. Then, we have a differential bigraded algebra $R = \Lambda[u_1, \ldots, u_m] \otimes A$ (Hereafter, $\otimes$ indicates $\otimes_\k$) with  map $d \colon R\to R$ by setting
$$
    \bideg u_i = (-1,2), \bideg v_i = (0,2), du_i =v_i, dv_i=0.
$$
Note that $R$ is a free $A$-module. Let $R^{-i} = \Lambda^{i}[u_1, \ldots, u_m] \otimes A$, where $\Lambda^{i}[u_1, \ldots, u_m]$ is the submodule of $\Lambda[u_1, \ldots, u_m]$ spanned by monomials of length $i$. Then, we have the following free resolution of $\k$, which is known as the \emph{Koszul resolution}:
$$ [R] \colon 0 \to R^{-m} \stackrel{d}{\to} \cdots \stackrel{d}{\to} R^{-1} \stackrel{d}{\to} A \stackrel{d}{\to} \k \to 0.$$

Let $P$ be an $n$-dimensional simple polytope with $m$ facets. The \emph{face ring} (or the \emph{Stanley-Reisner ring}) of $P$ is the quotient ring
$$
    \k(P) = \k[v_1, \ldots, v_m] / I_P,
$$
where $I_P$ is the homogeneous ideal generated by all square-free monomials $v_{i_1} v_{i_2} \cdots v_{i_s}$ such that $F_{i_1} \cap \cdots F_{i_s} = \emptyset$. The ideal $I_P$ is called the \emph{Stanley-Reisner ideal} of $P$. By identifying the polynomial ring $\k[v_1, \ldots, v_m]$ in the definition of $\k(P)$ with $A$ above, $\k(P)$ can be regarded as an $A$-module. By applying the functor $\otimes_A \k(P)$ to the Koszul resolution, we obtain the following cochain complex of graded modules:
$$
  [R\otimes_A \k(P)] \colon 0 \to R^{-m}\otimes_A \k(P) \to \cdots \to R^{-1}\otimes_A \k(P) \to \k(P),
$$
where the differential map is $d\otimes_A 1$. The $(-i)$-th cohomology module of the above cochain complex is denoted by $\Tor_A^{-i}(\k(P),\k)$, and we have the graded $A$-module $$\Tor_A(\k(P),\k) = \bigoplus_i \Tor_A^{-i}(\k(P),\k).$$
Note that there is a canonical multiplicative structure on
$$
    \Tor_A(\k(P),\k) = H[R \otimes_A \k(P)] = H[\Lambda[u_1, \ldots, u_m] \otimes \k(P)],
$$ and hence, $\Tor_A(\k(P),\k)$ is canonically a bigraded $\k$-algebra.
The bigraded algebra $\Tor^{\ast,\ast}_A(\k(P),\k)$ is called the \emph{Tor-algebra} of a simple polytope $P$, and the \emph{bigraded Betti numbers} of $P$ are defined by
$$
    \beta^{-i,2j}(P; \k) = \dim_\k \Tor^{-i,2j}_A(\k(P),\k) \quad \text{for $0 \leq i,j \leq m$.}
$$

Hereafter, only the case of $\k = \Q$ is considered. For simplicity, we set $\beta^{-i, 2j}(P) = \beta^{-i, 2j}(P;\Q)$.

The following theorem of Hochster \cite{Ho} gives a nice combinatorial interpretation of bigraded Betti numbers.
\begin{theorem}
Let $P$ be a simple convex polytope with facets $F_1, \ldots, F_m$. For a subset $\sigma \subset \{1, \ldots, m\}$, let $P_\sigma = \bigcup_{i\in\sigma}F_i \subset P$. Then, we have
$$
    \beta^{-i,2j}(P) = \sum_{|\sigma|=j} \dim \tilde{H}^{j-i-1}(P_\sigma;\Q).
$$ Here, $\dim \tilde{H}^{-1}(\emptyset)=1$ by convention.
\end{theorem}

\begin{example}
    Let $P$ be a $3$-dimensional simple polytope with $m$ facets $F_1, \ldots, F_m$.
    \begin{enumerate}
        \item $\beta^{-1,4}(P) = \sum_{1\leq i<j \leq m} \dim \tilde{H}^{0}(F_i \cup F_j;\Q)=$ the number of pairs of facets that do not intersect.
        \item $\beta^{-1,6}(P) = \sum_{1\leq i<j<k \leq m} \dim \tilde{H}^{1}(F_i \cup F_j \cup F_k;\Q)=$ the number of triple of facets whose union is homotopy equivalent to $S^1$. Such triple of facets is called a \emph{$3$-belt}.
        \item Assume $\beta^{-1,6}(P) = 0$.\footnote{A $3$-dimensional simple polytope $P$ whose $\beta^{-1,6}(P)$ is $0$ is said to be \emph{irreducible} because it cannot be expressed as a connected sum of a finite number of simple polytopes (see \cite{CK11}).} Then, since there is no $3$-belt, $\beta^{-2,8}(P)$ is equal to the number of quadruples of facets whose union is homotopy equivalent to $S^1$. Such quadruple of facets is called a \emph{$4$-belt}.
    \end{enumerate}
\end{example}

Let $P$ be a $3$-dimensional simple polytope with $m$ facets $F_1, \ldots, F_m$.
Now, we consider $\Lambda[u_1, \ldots, u_m] \otimes \k(P)$. Let $d$ be a differential operator on $\Lambda[u_1, \ldots, u_m] \otimes \Q(P)$ induced from $d\otimes_A 1$ on $[R \otimes_A \Q(P)]$. As mentioned before,
$\Tor^{\ast,\ast}_A(\Q(P),\Q) = H[\Lambda[u_1, \ldots, u_m] \otimes \k(P)]$.

Assume that $F_i$ and $F_j$ do not intersect. Then, $u_i v_j$ is an element of bidegree $(-1,4)$ in $\Lambda[u_1, \ldots, u_m] \otimes \Q(P)$. Since $d(u_i v_j) = v_i v_j = 0 \in \Q(P)$, it is a cycle. Furthermore, $u_i v_j$ and $v_i u_j$ are homologous because $d(u_iu_j)=v_iu_j - u_iv_j$, and $u_i v_j$ and $u_{i'}v_{j'}$ are not homologous for $\{i,j\} \neq \{i',j'\}$. Therefore, the set of equivalent classes $X^{-1,4}:=\{[u_iv_j] \mid F_i \cap F_j = \emptyset\}$ becomes a subset of generators of $H^{-1,4}[\Lambda[u_1, \ldots, u_m] \otimes \Q(P)]$. Since $|X^{-1,4}| = \beta^{-1,4}(P)$, the set $X^{-1,4}$ itself is the set of generators.

Assume that $\beta^{-1,6}(P)=0$ and $\{F_i, F_j, F_k, F_\ell\}$ is a $4$-belt of $P$. Then, $u_i v_j u_k v_\ell$ is an element of bidegree $(-2,8)$ in $\Lambda[u_1, \ldots, u_m] \otimes \Q(P)$. It may be assumed that $F_i \cap F_k = \emptyset$ and $F_j \cap F_\ell = \emptyset$. Therefore, $d (u_i v_j u_k v_\ell) = 0$. One can easily see that $u_i v_j u_k v_\ell$ is homologous to $u_{i'} v_{j'} u_{k'} v_{\ell'}$ if and only if $\{\{i',j'\},\{k',\ell'\}\}$ is equal to either $\{\{i,j\},\{k,\ell\}\}$ or $\{\{i,\ell\},\{k,j\}\}$ as a set of sets. Let $X^{-2,8}:=\{[u_i v_j u_k v_{\ell}] \mid \text{$\{F_i, F_j, F_k, F_\ell\}$ is a $4$-belt}\}$, which is a subset of generators of $H^{-2,8}[\Lambda[u_1, \ldots, u_m] \otimes \Q(P)]$. Since $|X^{-2,8}| = \beta^{-2,8}(P)$, the set $X^{-1,4}$ itself is the set of generators. Hence, we have the following proposition.

\begin{proposition}
Let $P$ be a $3$-dimensional simple polytope with $m$ facets $F_1, \ldots, F_m$, with $\beta^{-1,6}(P) =0$. Then, each generator of $\Tor_A^{-1,4}(\Q(P),\Q)$ can be indexed by a pair of facets that do not intersect, and each generator of $\Tor_A^{-2,8}(\Q(P),\Q)$ can be indexed by $4$-belts. Let $\sigma$ and $\tau$ be generators in $\Tor_A^{-1,4}(\Q(P),\Q)$ indexed by $\{F_1, F_2\}$ and $\{F_3, F_4\}$, respectively. Then,
  $$
        \sigma \cdot \tau = \left\{
                              \begin{array}{ll}
                                \eta \neq 0 & \hbox{if $\{F_1, F_2, F_3, F_4\}$ is a $4$-belt in $P$;} \\
                                0 & \hbox{otherwise,}
                              \end{array}
                            \right.
  $$ where $\eta$ is a generator in $\Tor_A^{-2,8}(\Q(P),\Q)$ indexed by $\{F_1, F_2, F_3, F_4\}$.
\end{proposition}

\section{Proof of Theorem~\ref{Theorem:Tor_is_not_decided_by_Betti_numbers}}\label{sec:proof_of_main_theorem}
It is to be noted that $\Tor^{\ast,\ast}_A(\Q(\cP),\Q)$ and $\Tor^{\ast,\ast}_A(\Q(\cQ),\Q)$ are isomorphic as groups. Hence, their multiplicative structures should be compared. We index each facet of $\cP$ and $\cQ$ as shown in Figure~\ref{Figure:indexing_of_facets_of_P_and_Q}.
\begin{figure}
\begin{pspicture}(0,0)(5,5)
\pspolygon(2.5,0.2)(0.5,1.1)(0.2,1.9)(0.5,3.2)(1.4,4.4)(2.5,4.8)(3.6,4.4)(4.5,3.2)(4.8,1.9)(4.5,1.1)
\psline(0.5,3.2)(2.5,2.5)\psline(4.5,3.2)(2.5,2.5)
\psline(0.2,1.9)(1.9,1.4)\psline(4.8,1.9)(3.1,1.4)
\psline(2.5,2.5)(2.5,2)\psline(2.5,0.2)(2.5,1)
\pspolygon(2.5,1)(1.9,1.4)(2.5,2)(3.1,1.4)
\psline[linestyle=dotted] (1.4,4.4)(1.4,3.1)
\psline[linestyle=dotted] (2.5,4.8)(2.5,3.5)
\psline[linestyle=dotted] (3.6,4.4)(3.6,3.1)
\psline[linestyle=dotted] (0.5,1.1)(1.4,3.1)(2.5,3.5)(3.6,3.1)(4.5,1.1)
\rput(1.8,4.7){$a$}
\rput(0.3,3){$b$}
\rput(2.5,3.3){$c$}
\rput(3.2,4.7){$d$}
\rput(3.5, 0.4){$e$}
\rput(1.5,1){$f$}
\rput(1.5,2){$g$}
\rput(3.5,2){$h$}
\rput(4.7,3){$i$}
\rput(3.5,1){$j$}
\rput(2.5,1.4){$k$}
\end{pspicture}
\begin{pspicture}(0,0)(5,5)
\pspolygon(3.6,0.4)(4.7,1.2)(4.9,2.5)(4.4,3.9)(3.5,4.7)(1.5,4.7)(0.6,3.9)(0.1,2.5)(0.3,1.2)(1.4,0.4)
\psline(4.4,3.9)(3.4,3) \psline(0.6,3.9)(1.6,3)
\psline(4.9,2.5)(3.6,1.6) \psline(0.1,2.5)(1.4,1.6)
\psline(3.6,0.4)(3.6,1.6)(3.1,2.2) \psline(1.4,0.4)(1.4,1.6)(1.9,2.2)
\pspolygon(3.1,2.2)(3.4,3)(1.6,3)(1.9,2.2)
\psline[linestyle=dotted] (4.7,1.2)(3.5,2.7)
\psline[linestyle=dotted] (0.3,1.2)(1.5,2.7)
\psline[linestyle=dotted] (3.5,4.7)(3.5,2.7)(1.5,2.7)(1.5,4.7)
\rput(4.6,3.9){$a$}
\rput(2.5,4.9){$b$}
\rput(2.5,3.9){$c$}
\rput(4,2.5){$d$}
\rput(4.1, 1.5){$e$}
\rput(2.5,0.2){$f$}
\rput(0.4,3.9){$g$}
\rput(1,2.5){$h$}
\rput(2.5,2.5){$i$}
\rput(2.5,1.3){$j$}
\rput(0.9,1.5){$k$}
\end{pspicture}
 \caption{Indices of facets of $\cP$ and $\cQ$}  \label{Figure:indexing_of_facets_of_P_and_Q}
\end{figure}

Let $A=\Q[a,b,\ldots, k]$, where the letters are degree $2$ indeterminates corresponding to the facets of $\cP$.
Now, consider the subspace
$$
    V_P:=\{ x \in \Tor_A^{-1,4}(\Q(P),\Q) \mid x r = 0 \text{ for all $r \in \Tor_A^{-1,4}(\Q(P),\Q)$}\}
$$ as a vector space over $\Q$. It is obvious that the dimension of $V_P$ is a ring invariant of $\Tor_A(\Q(P),\Q)$.

Note that $\beta^{-2,8}(\cP)=\beta^{-2,8}(\cQ)=4$. In other words, there are four $4$-belts in both $\cP$ and $\cQ$.
In $\cP$, all $4$-belts are indexed by $\{b,c, d,e\}$, $\{g,f,j,h\}$, $\{a,c,i,e\}$, and $\{c,i,e,b\}$.
Hence, only some of the products between two pairs among $\{b,d\}$, $\{c,e\}$, $\{g,j\}$, $\{f,h\}$, $\{a,i\}$, and $\{b,i\}$ ($6$ generators) can be nonzero. This implies that $\dim_\Q V_P = 28 - 6 = 22$.

In $\cQ$, all $4$-belts are indexed by $\{a,c,g,f\}$, $\{a,d,j,f\}$, $\{c,h,j,d\}$, and $\{f,g,h,j\}$.
Hence, only some of the products between two pairs among $\{a,g\}$, $\{c,f\}$, $\{a,j\}$, $\{d,f\}$, $\{c,j\}$, $\{h,d\}$, $\{f,h\}$, and $\{g,j\}$ ($8$ generators) can be nonzero. This implies that $\dim_\Q V_Q = 28 - 8 = 20$. Therefore, $\Tor_A(\Q(\cP),\Q)$ and $\Tor_A(\Q(\cQ),\Q)$ are not isomorphic as rings, which proves the first part of the theorem.

We note that $H^3(\cZ_P;\Q) = \Tor_A^{-1,4}(\Q(P),\Q)$ and $$H^6(\cZ_P;\Q) = \Tor_A^{-2,8}(\Q(P),\Q) \oplus \Tor_A^{-4,10}(\Q(P),\Q).$$
Since any element of $\Tor_A^{-4,10}(\Q(P),\Q)$ cannot be expressed as two elements in $\Tor_A^{-1,4}(\Q(P),\Q)$ because of degrees, one can show that $H^\ast(\cZ_\cP;\Q) \not \cong H^\ast(\cZ_\cQ;\Q)$ as rings by using the above argument, which proves the second part of the theorem.

\section{Proof of Theorem~\ref{Theorem:cohom_does_not_imply_combi}}
First, it is obvious that $\cP$ and $\cQ$ are not combinatorially rigid since they have the same bigraded Betti numbers while they are not combinatorially equivalent.

In the remaining part of this section, we prove that $\cP$ and $\cQ$ are cohomologically rigid.

Suppose that there exist a polytope $R$ and quasitoric manifolds $M$ and $N$ over $\cP$ and $P'$, respectively, such that $H^\ast(M)$ and $H^\ast(N)$ are isomorphic as graded rings. Then, by \cite[Lemma 3.7]{ch-pa-su10}, it follows that $P'$ has $11$ facets, and $$\Tor^{\ast,\ast}_A(\Q(\cP),\Q) = \Tor^{\ast,\ast}_A(\Q(P'),\Q).$$
In particular, $\beta^{\ast,\ast}(\cP) = \beta^{\ast,\ast}(P')$.

Now, let us investigate all other polytopes with $11$ facets. A graph $G$ is said to be \emph{$P^3$-realizable} if there is a $3$-dimensional polytope whose corresponding $1$-complex is isomorphic to $G$. Let $P(G)$ denote such a polytope. A graph $G$ is said to be \emph{$k$-(vertex-)connected} if there is no set of $k-1$ vertices that, when removed, disconnects the graph. It is known that a graph $G$ is planar and $3$-connected if and only if $G$ is $P^3$-realizable (\cite{St}). A $P^3$-realizable graph is called a \emph{triangulation} if all the faces of the graph are triangle when the graph is embedded into a $2$-dimensional sphere $S^2$. Hence if $G$ is a $3$-connected triangulation, then $P(G)$ is a simplicial polytope that is dual to a simple polytope.

Using the graph-generating program \emph{plantri} developed by Brinkmann and McKay, we can list all $3$-connected triangulations with a certain number of vertices. Such a list gives us all $3$-dimensional simple polytopes with a certain number of facets. Using the program \emph{Macaulay2} again, we can list all $3$-dimensional simple polytopes $P$ with $11$ facets satisfying $\beta^{-1,6}(P)=0$ and compute their bigraded Betti numbers (see \cite{Y.Ch}). See Table~\ref{table:11facets}; each polytope has $11$ facets $a, b, \ldots, k$. Each polytope is indexed by using information of adjacency of the facets. The $n$-th component is the list of facets that intersect the (alphabetical) $n$-th facet. The Betti numbers are listed in the form
$$
    (\beta^{-1,4}(P), \ldots, \beta^{-(j-1), 2j}(P), \ldots, \beta^{-7,16}(P)).
$$
Note that the integer tuple of the above form completely determines all the bigraded Betti numbers of a $3$-dimensional polytope (see \cite[Section~7]{ch-pa-su10} for details).

\begin{table}
  {\Tiny
    \centering
\begin{tabular}{|c|c|c|}
  \hline
   & Polytope & Betti number \\ \hline

1 & bcdef,afghc,abhid,acijke,adkf,aekgb,bfkjh,bgjic,chjd,dihgk,djgfe &(28, 105, 164, 112, 28, 2, 0) \\
2 & bcdef,afghijc,abjd,acjke,adkgf,aegb,bfekh,bgki,bhkj,bikdc,djihge&(28, 105, 167, 131, 47, 5, 0) \\
3 & bcde,aefghic,abid,acijke,adkfb,bekjg,bfjh,bgji,bhjdc,dihgfk,djfe&(28, 105, 169, 138, 54, 7, 0) \\
4 & bcde,aefghijc,abjd,acjke,adkfb,bekg,bfkh,bgki,bhkj,bikdc,djihgfe&(28, 105, 175, 159, 75, 13, 0) \\
5 & bcdef,afghijc,abjd,acje,adjkf,aekgb,bfkh,bgki,bhkj,bikedc,ejihgf&(28, 105, 172, 144, 60, 10, 0) \\
6 & bcde,aefc,abfgd,acghe,adhijfb,bejgc,cfjkhd,dgkie,ehkj,eikgf,gjih&(28, 105, 171, 141, 57, 9, 0) \\
7 & bcde,aefgc,abghijkd,acke,adkjfb,bejihg,bfhc,cgfi,chfj,cifek,cjed&(28, 105, 174, 156, 72, 12, 0) \\
8 & bcde,aefghc,abhijd,acje,adjfb,bejkg,bfkh,bgkic,chkj,cikfed,fjihg&(28, 105, 168, 129, 45, 6, 0) \\
9 & bcde,aefghc,abhijd,acje,adjfb,bejikg,bfkh,bgkic,chkfj,cifed,fihg&(28, 105, 170, 136, 52, 8, 0) \\
10 & bcdef,afghic,abid,acijke,adkf,aekgb,bfkjh,bgji,bhjdc,dihgk,djgfe&(28, 105, 165, 119, 35, 3, 0) \\
11 & bcde,aefghic,abid,acije,adjkfb,bekg,bfkh,bgkji,bhjdc,dihke,ejhgf&(28, 105, 170, 136, 52, 8, 0) \\
12 & bcde,aefgc,abghid,acie,adijfb,bejkg,bfkhc,cgkji,chjed,eihkf,fjhg&(28, 105, 166, 123, 39, 4, 0) \\
13 & bcde,aefgc,abghid,acie,adijfb,bejkg,bfkhc,cgki,chkjed,eikf,fjihg&(28, 105, 167, 125, 41, 5, 0) \\
14 & bcde,aefghc,abhd,achije,adjfb,bejikg,bfkh,bgkidc,dhkfj,dife,fihg&(28, 105, 169, 134, 50, 7, 0) \\
15 & bcde,aefghc,abhijd,acjgfe,adfb,bedg,bfdjkh,bgkic,chkj,cikgd,gjih&(28, 105, 173, 145, 61, 11, 0) \\
16 & bcde,aefc,abfghid,acie,adijkfb,bekgc,cfkjh,cgji,chjed,eihgk,ejgf&(28, 105, 170, 143, 59, 8, 0) \\
17 & bcde,aefc,abfghid,acie,adihjfb,bejkgc,cfkh,cgkjei,ched,ehkf,fjhg&(28, 105, 177, 159, 75, 15, 0) \\
18 & bcde,aefghic,abid,acijgke,adkfb,bekg,bfkdjh,bgji,bhjdc,dihg,dgfe&(28, 105, 173, 149, 65, 11, 0) \\
19 & bcde,aefghijc,abjd,acjkhgfe,adfb,bedg,bfdh,bgdki,bhkj,bikdc,djih&(28, 105, 179, 169, 85, 17, 0) \\
20 & bcde,aefghijkc,abkd,ackjihgfe,adfb,bedg,bfdh,bgdi,bhdj,bidk,bjdc&(28, 105, 189, 189, 105, 27, 0) \\
21 & bcde,aefgc,abghd,ache,adhijfb,bejg,bfjkhc,cgkied,ehkj,eikgf,gjih&(28, 105, 171, 141, 57, 9, 0) \\
22 & bcde,aefc,abfghd,ache,adhijfb,bejkgc,cfkh,cgkied,ehkj,eikf,fjihg&(28, 105, 173, 145, 61, 11, 0) \\
23 & bcdefg,aghc,abhijd,acje,adjf,aejkhg,afhb,bgfkic,chkj,cikfed,fjih&(28, 105, 171, 137, 53, 9, 0) \\
24 & bcdef,afgc,abghid,acije,adjf,aejkgb,bfkhc,cgkji,chjd,dihkfe,fjhg&(28, 105, 166, 123, 39, 4, 0) \\
25 & bcdef,afgc,abghijd,acje,adjkhgf,aegb,bfehc,cgeki,chkj,ciked,ejih&(28, 105, 173, 149, 65, 11, 0) \\
  \hline
\end{tabular}
}
  \caption{Irreducible polytopes with 11 facets}\label{table:11facets}
\end{table}

In Table~\ref{table:11facets}, the $11$-th polytope is $\cP$, and the $24$-th polytope is $\cQ$. One can easily check that there is no other polytope whose bigraded Betti numbers are equal to those of $\cP$ and $\cQ$. Thus, $P'$ cannot be combinatorially equivalent to any polytope with $11$ facets other than $\cP$ and $\cQ$. Moreover, by Theorem~\ref{Theorem:Tor_is_not_decided_by_Betti_numbers}, $P'$ cannot be combinatorially equivalent to $\cQ$.
Therefore, $P'$ is $\cP$, which proves that $\cP$ is cohomologically rigid.

Similar arguments can be presented for $\cQ$ to prove its cohomological rigidity.

\bigskip
\bibliographystyle{amsalpha}

\begin{thebibliography}{CPS08}

\bibitem{BP}
V.~M. Buchstaber and T.~E. Panov, \emph{Torus actions and their
  applications in topology and combinatorics}, University Lecture Series,
  vol.~24, American Mathematical Society, Providence, RI, 2002.

\bibitem{CK11}
S. Choi and J.~S. Kim, \emph{Combinatorial Rigidity of 3-dimensional simplicial polytopes},  Int. Math. Res. Not. IMRN. 2011(8), (2011), 1935--1951.

\bibitem{ch-ma-su11}
S. Choi, M. Masuda and D.~Y. Suh, \emph{Rigidity problems in toric topology, a survey}, Proc. Steklov Inst. Math., 275 (2011), 177--190.

\bibitem{ch-pa-su10}
S. Choi, T.~E. Panov, and D.~Y. Suh, \emph{Toric cohomological
  rigidity of simple convex polytopes},  J. London Math. Soc. 82(2) (2010), 343--360.

\bibitem{Y.Ch}
Y. Choi, \emph{Cohomological rigidity of simple $3$-polytopes with $10$ facets}, Master Thesis, KAIST, (2008).

\bibitem{DJ}
M.l~W. Davis and T. Januszkiewicz, \emph{Convex polytopes, {C}oxeter
  orbifolds and torus actions}, Duke Math. J. \textbf{62} (1991), no.~2,
  417--451.

\bibitem{Ho}
M. Hochster, \emph{Cohen-Macaulay rings, combinatorics, and simplicial complexes}, in: Ring Theory II (Proc. Second Oklahoma Conference), B. R. McDonald and R. Morris, eds.,Dekker, New York, (1977), 171--223.

\bibitem{St}
E. Steinitz, \emph{Polyeder und Raumeinteilungen}, Enzykl, Math. Wiss, Vol. 3 (Geometrie)
Part 3AB12, (1922), 1–-139.
\end{thebibliography}

\end{document}